\documentclass[12pt,a4paper]{amsart}
\usepackage[a4paper,left=28mm,right=28mm,top=28mm,bottom=34mm]{geometry}
\usepackage{amsmath,amssymb,amsthm,mathtools}
\usepackage[T1]{fontenc}
\usepackage{lmodern}
\usepackage[colorlinks=true,linkcolor=blue,citecolor=blue,urlcolor=blue]{hyperref}

\newtheorem{theorem}{Theorem}[section]
\newtheorem{lemma}[theorem]{Lemma}
\newtheorem{proposition}[theorem]{Proposition}
\newtheorem{corollary}[theorem]{Corollary}
\theoremstyle{remark}
\newtheorem{remark}[theorem]{Remark}

\newcommand{\F}{\mathbb F}
\newcommand{\Nrd}{\operatorname{Nrd}}
\newcommand{\Nm}{\operatorname{N}}
\newcommand{\Tr}{\operatorname{Tr}}
\newcommand{\SL}{\operatorname{SL}}
\newcommand{\GL}{\operatorname{GL}}
\newcommand{\rad}{\operatorname{rad}}
\newcommand{\ord}{\operatorname{ord}}
\newcommand{\Gal}{\operatorname{Gal}}

\title[Counting irreducible skew polynomials]{%
Counting irreducible skew polynomials\\$x^{2n}+ax^n+1$}

\author[Tomoya Hayata]{Tomoya Hayata}
\address[Tomoya Hayata]{%
Keio University School of Medicine,
4-1-1 Hiyoshi, Kohoku-ku, Yokohama, 223-8521, Japan}
\address[Tomoya Hayata]{%
RIKEN Center for Interdisciplinary Theoretical
and Mathematical Sciences (iTHEMS), RIKEN,
Wako 351-0198, Japan}
\email{hayata@keio.jp}

\author[Takahiro Tsushima]{Takahiro Tsushima}
\address[Takahiro Tsushima]{Keio University School of Medicine,
4-1-1 Hiyoshi, Kohoku-ku, Yokohama, 223-8521, Japan}
\email{tsushima@keio.jp}

\makeatletter
\@ifundefined{subjclassname@2020}{%
  \@namedef{subjclassname@2020}{\textup{2020} Mathematics Subject Classification}%
}{}
\makeatother


\thanks{The first author was supported by JSPS KAKENHI Grant Number
JP24K00630 and JST COI-NEXT Grant Number JPMJPF2221.
The second author was supported by JSPS KAKENHI Grant Number JP25K06959.}
\date{}

\hypersetup{
  pdftitle={Counting irreducible skew polynomials},
  pdfauthor={Tomoya Hayata and Takahiro Tsushima}
}

\begin{document}

\begin{abstract}
We give a counting formula for the coefficients $a\in\F_{p^h}$ for which
$x^{2n}+ax^n+1$ is irreducible in $\F_{p^h}[x;\sigma]$, for arbitrary positive integers $n$ and $h$, where
$\sigma(a)=a^p$ and $p$ is an arbitrary prime power.
The proof reduces the problem to quadratic reduced norms, establishes a
uniformity statement for every Frobenius generator, and counts the
remaining coefficients by M\"obius inversion.
\end{abstract}

\maketitle

\section{Introduction}\label{sec:introduction}

Throughout, $p$ denotes a prime power, not necessarily the characteristic.
Let $q=p^h$, and let
$\sigma\in\operatorname{Aut}(\F_q)$ be given by $\sigma(c)=c^p$.
We study the number of coefficients $a\in\F_q$ for which
\[
 f_a(x)=x^{2n}+ax^n+1
\]
is irreducible in the skew polynomial ring $\F_q[x;\sigma]$, whose
multiplication is determined by $xc=\sigma(c)x$.
Our aim is an explicit formula valid for every pair of positive
integers $n,h$, including the case $n\nmid h$.

The main structural result is a uniformity theorem for quadratic
reduced norms on the one-parameter family $t^2+at+1$.
Let $P$ be any prime power, and let $\tau(c)=c^{P^r}$ act on
$\F_{P^k}$, with $\gcd(r,k)=1$.
If
\[
 \Nrd_{\F_{P^k}[t;\tau]}(t^2+at+1)
       =Z^2+b_r(a)Z+1,\qquad Z=t^k,
\]
then for every $b\in\F_P$ such that $Z^2+bZ+1$ is irreducible,
\[
 \#\{a\in\F_{P^k}:b_r(a)=b\}=D_k(P).
\]
Here $D_k(P)$ is the explicit quantity defined in \eqref{eq:D}.
Thus the norm map has equally sized fibers over all irreducible
quadratics with constant term one, and the same fiber size holds
for every generator of $\Gal(\F_{P^k}/\F_P)$.
This separates the count of the original coefficients from the
arithmetic restrictions on their norm values.

To state the resulting counting formula, put
\[
 d=\gcd(n,h),\qquad m=n/d,\qquad H=h/d,\qquad P=p^d.
\]
Writing $N_n(p,h)$ for the number of irreducible members of the family,
we prove
\begin{equation}\label{eq:intro-main}
 N_n(p,h)=D_H(P)B_{d,m}(p).
\end{equation}
The factor $D_H(P)$ is given by
\[
 D_H(P)=
 \begin{cases}
 (P^H+1)/(P+1),&H\text{ odd},\\[2pt]
 \bigl(1-(-P)^{H/2}\bigr)^2/(P+1),&H\text{ even}.
 \end{cases}
\]
The factor $B_{d,m}(p)$ counts coefficients $b\in\F_{p^d}$ for which
$V^2+bV+1$ is irreducible, $b$ has degree $d$ over $\F_p$, and
$X^m-\alpha$ is irreducible over $\F_{p^{2d}}$ for a root $\alpha$ of
the quadratic. Theorem~\ref{thm:main} evaluates this factor by an
explicit divisor sum. Thus \eqref{eq:intro-main} can be evaluated using integer
arithmetic, without enumerating the elements of $\F_q$.

For the standard generator, we express the fiber as a twisted
Frobenius fixed-point count on the Drinfeld curve
\[
 XY^P-X^PY=1.
\]
For background on this curve and its group action, see
\cite[Chapter~2]{Bonnafe}; the geometric facts and trace calculations
needed here are verified in Section~\ref{sec:standard}.
For a general generator, suitable normal forms for two-dimensional
semilinear maps reduce the comparison to power maps having the same
image and the same fiber size.
Combining this uniformity with the norm reduction and M\"obius inversion
gives \eqref{eq:intro-main}.

The distinction from root-count distributions is important here.
For the standard Frobenius, the linearized polynomial
$X^{P^2}+aX^P+X$ leads, by $Y=X^{P-1}$ for $X\ne0$, to the projective
polynomial $Y^{P+1}+aY+1$ studied in Bluher's framework~\cite{Bluher}.
The number of $K$-rational roots of this projective polynomial is the
number of $\F_P$-lines in the root space invariant under $F^k$,
where $K=\F_{P^k}$ and $F$ is the $P$-power Frobenius.
This number does not, in general, determine the characteristic
polynomial of $F^k$. Thus an aggregate root-count distribution does not by
itself give the count for each prescribed irreducible quadratic norm.
Our curve argument keeps track of that characteristic polynomial.

Skew polynomial rings were introduced by Ore~\cite{Ore} and provide
an algebraic setting for division, factorization, and semilinear
operators. For finite coefficient fields, Caruso and Le Borgne identify
reduced norms with characteristic polynomials of iterated companion
semilinear maps and characterize irreducibility through the norm
\cite[Lemma 2.1.15 and Proposition 2.1.17]{CL}.
Pumpl\"un and Thompson~\cite{PT} study reduced norms and their relation
to reducibility for a broader class of skew polynomial rings.
For more general Ore extensions, Churchill and Zhang~\cite{CZ} obtain
irreducibility criteria using noncommutative valuations.
In the finite-field setting considered here, the reduced-norm criterion
provides the structural reduction. To obtain an exact coefficient
count, one must additionally determine the fibers of the norm on the
prescribed family and impose the conditions arising from power
substitution. It is this restricted fiber count that is supplied by
the uniformity theorem above.

Odoni~\cite[Theorem 1]{Odoni} counts all irreducible additive
polynomials of a given degree under composition, using the corresponding
skew polynomial ring and its central quotients. Such a total count
does not determine how many irreducible elements have specified
coefficients or a prescribed sparse form. Our family fixes the leading
and constant coefficients at one and allows only the term $ax^n$
between them. The passage from total enumeration to this constrained
count requires the norm fibers on that family.
Factorization algorithms over finite fields were developed by
Giesbrecht~\cite{Giesbrecht} and improved by Caruso and Le Borgne~\cite{CL}.
Skew polynomial rings also occur in the construction of skew-cyclic
codes~\cite{BGU}; a survey of the ring-theoretic and coding-theoretic
background is given in \cite{GLSurvey}.

A motivation for the prescribed family comes from the second author's study of self-adjoint
irreducible polynomials in skew Laurent polynomial rings over finite
fields~\cite{Tsushima}.
That work determines the complete fiber-size distribution of
$c\mapsto c^{p^j}+c^{-1}$ for arbitrary positive integers $j$
\cite[Corollary 3.12]{Tsushima}, by combining Bluher's
root-counting results~\cite{Bluher} with a second-moment argument
based on point counts on Edixhoven curves.
In particular, \cite[Corollary 5.11]{Tsushima} gives the number of
monic self-adjoint irreducible Laurent polynomials of degree two,
equivalently the total number of irreducible polynomials $x^2+ax+1$,
with $D_h(p)$ as in \eqref{eq:D}:
\[
 N_1(p,h)=\frac{p-\epsilon_p}{2}D_h(p),
 \qquad
 \epsilon_p=\begin{cases}0,&p\text{ even},\\1,&p\text{ odd}.\end{cases}
\]
The same work also gives the total degree-two count for arbitrary
Frobenius powers \cite[Corollary 4.3]{Tsushima}.
The present paper refines these total counts by prescribing each
irreducible quadratic reduced norm and proving independence of the
Frobenius generator for the individual fibers. These refined counts
allow us to impose the exact-degree and power-substitution conditions
needed to treat $x^{2n}+ax^n+1$ for all $n,h$.
For this same sparse family, \cite[Corollary 4.2]{Tsushima} gives an
upper bound by excluding coefficients that admit a factorization
$(x^n+b)(x^n+c)$. Here we obtain the exact count by imposing
necessary and sufficient conditions for irreducibility.
The general higher-degree bound in \cite[Proposition 6.8]{Tsushima}
concerns the larger family of self-adjoint Laurent polynomials.

Indeed, multiplication on the left by $x^{-n}$ gives
\[
 x^{-n}f_a(x)=x^n+\sigma^{-n}(a)+x^{-n}.
\]
Since $a\mapsto\sigma^{-n}(a)$ is a bijection of $\F_q$, and passing
from the polynomial ring to the Laurent ring preserves irreducibility
for polynomials with nonzero constant term, our count is equivalently
an exact count of the self-adjoint irreducible Laurent polynomials
$x^n+c+x^{-n}$, with $c\in\F_q$.

A further motivation for counting self-adjoint irreducible Laurent
polynomials comes from their connection with primitive local Galois
representations. In the construction of \cite{T}, irreducibility of
an associated additive polynomial under composition gives a sufficient
condition for the resulting representation to be primitive
\cite[Corollary 3.28(3)]{T}.

The arithmetic part is related to the classical enumeration of
self-reciprocal irreducible polynomials over finite fields.
Meyn and G\"otz~\cite{MG} give a M\"obius-inversion proof of that
enumeration and describe its relation to the groups of roots of unity
of order $p^j+1$.
In our setting, a root $\alpha$ of an irreducible quadratic
$V^2+bV+1\in\F_{p^d}[V]$ satisfies
\[
 \alpha^{p^d+1}=1,
 \qquad b=-(\alpha+\alpha^{-1}).
\]
This connects exact-degree conditions on the coefficient $b$ with
self-reciprocal polynomials. The remaining power-substitution condition
is expressed by the usual finite-field irreducibility criterion for
$X^m-\alpha$; see, for example, \cite[Theorem 3.75]{LN}.

When $h=nk$, Corollary~\ref{cor:multiple} gives
\[
 N_n(p,nk)=C_n(p)D_k(p^n)
         =nS_p(2n)D_k(p^n),
\]
where $C_n(p)=B_{n,1}(p)$ and $S_p(2n)$ denotes the number of monic
self-reciprocal irreducible polynomials of degree $2n$ over $\F_p$.
The identity $C_n(p)=nS_p(2n)$ is explained in
Section~\ref{sec:specializations}; see \eqref{eq:selfreciprocal-count}.
This identifies the classical self-reciprocal count inside the formula
and isolates the extension factor $D_k(p^n)$.
The specializations $n=2,3,4$ are made explicit at the end of the paper.
If $n$ is a power of two, the formula also gives $N_n(p,h)=0$ whenever $n\nmid h$.

Generative AI assistance was used in formulating the conjecture for
general $n$, exploring and developing proof arguments, and preparing
the exposition.
It was also used to develop code for supplementary numerical
checks, which are not used in the proofs.
Further details are provided in the declaration at the end of the paper.

The paper is organized as follows.
Section~\ref{sec:statement} states the arithmetic formula.
Section~\ref{sec:reduction} reduces irreducibility to a quadratic norm
and a binomial condition.
Sections~\ref{sec:standard} and~\ref{sec:generator} establish the fiber
count for the standard and arbitrary Frobenius generators,
respectively.
Section~\ref{sec:arithmetic} completes the arithmetic count and proves
the main theorem. Section~\ref{sec:specializations} records its
specializations.

\section{Statement of the formula}\label{sec:statement}
Let $p$ be a prime power, $n,h\geq1$, and put
\[
 N_n(p,h)=\#\left\{a\in\F_{p^h}:x^{2n}+ax^n+1
       \text{ is irreducible in }\F_{p^h}[x;\sigma]\right\},
 \qquad \sigma(c)=c^p.
\]
Set
\[
 d=\gcd(n,h),\qquad m=n/d,\qquad H=h/d,\qquad P=p^d.
\]
For any prime power $P$ define
\begin{equation}\label{eq:D}
 D_k(P)=
 \begin{cases}
 (P^k+1)/(P+1),&k\text{ odd},\\[2pt]
 \bigl(1-(-P)^{k/2}\bigr)^2/(P+1),&k\text{ even}.
 \end{cases}
\end{equation}
Write $R=\rad(m)$, with $\rad(1)=1$. If $m$ is even or $R\nmid P+1$,
define $B_{d,m}(p)=0$. Otherwise set
\begin{equation}\label{eq:B}
 L=\prod_{\ell\mid R}\ell^{v_\ell(P+1)},\qquad
 g_e=\gcd(P+1,p^e-1),\qquad
 B_{d,m}(p)=\frac{\varphi(R)}{2R}
   \sum_{\substack{e\mid2d\\L\mid g_e}}
        \mu(2d/e)g_e.
\end{equation}
Empty products are $1$; in particular $R=L=1$ when $m=1$.
Here $\mu$ and $\varphi$ denote the M\"obius and Euler functions,
$\rad(m)$ is the product of the distinct prime divisors of $m$,
and $v_\ell$ denotes the $\ell$-adic valuation.

\begin{theorem}\label{thm:main}
With the preceding notation,
\[
 N_n(p,h)=D_H(p^d)\,B_{d,m}(p).
\]
\end{theorem}

The two factors in the formula have different roles. The factor
$B_{d,m}(p)$ counts the quadratic norm coefficients $b\in\F_P$ that
satisfy the irreducibility conditions obtained in
Section~\ref{sec:reduction}. For each such $b$, the factor $D_H(P)$
counts the coefficients $a\in\F_{p^h}$ having that norm coefficient.
Thus the proof separates into an arithmetic count of admissible norm
coefficients and a uniform count of their preimages. The latter is
established in Sections~\ref{sec:standard} and~\ref{sec:generator},
and the two counts are combined in Section~\ref{sec:arithmetic}.

We fix our convention for reduced norms.
Let $K$ be a finite field, let $\tau\in\operatorname{Aut}(K)$
have order $s$, and put
\[
 K_0=K^\tau,\qquad \mathcal R=K[t;\tau],\qquad z=t^s.
\]
The center of $\mathcal R$ is $K_0[z]$, and its central localization
\[
 \mathcal A
 =\mathcal R\otimes_{K_0[z]}K_0(z)
\]
is a central simple algebra of degree $s$ over $K_0(z)$.
For $f\in\mathcal R$, we write
\[
 \Nrd_{\mathcal R}(f)
 :=\Nrd_{\mathcal A/K_0(z)}(f)\in K_0[z].
\]
Here the reduced norm is defined by extending scalars to a
splitting field, identifying the resulting algebra with a
matrix algebra, and taking the determinant.

Concretely, $\mathcal A$ is a left $K(z)$-vector space
with basis $\{1,t,\ldots,t^{s-1}\}$.
Right multiplication by $f$ is $K(z)$-linear, and
$\Nrd_{\mathcal R}(f)$ is the determinant of its matrix in this basis.
We always specify the central variable when comparing
reduced norms in different skew polynomial rings.

We record explicitly how the determinant convention relates to
companion semilinear maps. This is the even-degree form of the standard
norm identity and irreducibility criterion; compare
\cite[Lemma 2.1.15 and Proposition 2.1.17]{CL}.
Let $f=t^D+\sum_{i<D}c_it^i$ be monic of even degree $D$ with $c_0\ne0$.
On the left $K$-vector space
\[
 M_f=\mathcal R/\mathcal Rf,
\]
left multiplication by $t$ induces a $\tau$-semilinear map $\Phi$.
In the basis $\{1,t,\ldots,t^{D-1}\}$ its coefficient matrix is the
companion matrix $C_f$, with last column
$(-c_0,\ldots,-c_{D-1})^t$. In particular,
\[
 \Phi^s=C_f\tau(C_f)\cdots\tau^{s-1}(C_f)
\]
as a $K$-linear map. With $z=t^s$, we have
\[
 \Nrd_{\mathcal R}(f)=\det(zI_D-\Phi^s),
 \qquad
 f\text{ irreducible}\ \Longleftrightarrow\
 \Nrd_{\mathcal R}(f)\text{ irreducible in }K_0[z].
\]

\begin{proof}[Verification of the conventions and criterion]
View $\mathcal R$ as a free left $K[z]$-module with basis
$\{1,t,\ldots,t^{s-1}\}$. Right multiplication by $f$ is injective and
has cokernel $M_f$. Smith normal form over $K[z]$ shows that its
determinant, up to a nonzero constant, is the characteristic polynomial
of multiplication by $z$ on $M_f$, namely that of $\Phi^s$.
The leading coefficient of this determinant is
$(-1)^{D(s-1)}$: its leading term comes from right multiplication by
$t^D$, and the matrix for right multiplication by $t$ has determinant
$(-1)^{s-1}z$. More explicitly, a determinant term formed from
monomials of degrees $j_0,\ldots,j_{s-1}\le D$ has $z$-degree
$(j_0+\cdots+j_{s-1})/s$; degree $D$ requires all $j_i=D$.
Since $D$ is even, the leading coefficient is one. This proves the
displayed identity, including the order of the twisted matrix product.

For completeness, the same determinant calculation gives
$\deg_z\Nrd(g)=\deg_tg$ for any nonzero skew polynomial $g$.
Multiplicativity therefore shows that an irreducible norm forces
$f$ to be irreducible. Conversely, if $f$ is irreducible, $M_f$ is
a simple $\mathcal R$-module. Its central annihilator is generated
by an irreducible polynomial $\eta\in K_0[z]$; this follows since
the image of the center in the endomorphism division ring of $M_f$
is a finite integral domain. As $c_0\ne0$, $\Phi$ is invertible,
so $\eta\ne z$.
Put $E=K_0[z]/(\eta)$ and $e=\deg\eta$.
The central quotient $\mathcal R/\eta\mathcal R$ is a central
simple algebra of degree $s$ over the finite field $E$, hence
is $M_s(E)$; see also \cite[Theorem 2.1.2 and Corollary 2.1.3]{CL}.
Its simple module has $E$-dimension $s$. Comparing $K_0$-dimensions
gives $sD=se$, hence $D=e$.
Since $\eta$ annihilates $\Phi^s$, and the characteristic polynomial
of $\Phi^s$ belongs to $K_0[z]$, all its irreducible factors over
$K_0$ equal $\eta$. Equality of degrees and monicity therefore give
$\Nrd(f)=\eta$.
\end{proof}

For $f=t^2+at+1$, this calculation uses exactly
\[
 C_f=A_a=\begin{pmatrix}0&-1\\1&-a\end{pmatrix},
 \qquad
 \Phi^s=A_a\tau(A_a)\cdots\tau^{s-1}(A_a).
\]
Each factor has determinant one.
All reduced norms below use the indicated central variable.

\section{Reduction to a quadratic norm}\label{sec:reduction}
The aim of this section is to replace irreducibility in the original
skew polynomial ring by three conditions over finite fields. We first
express the reduced norm in terms of a quadratic polynomial with
coefficient $b_a\in\F_P$. We then show that irreducibility is equivalent
to the irreducibility of this quadratic, an exact-degree condition on
$b_a$, and the irreducibility of a binomial over $\F_{P^2}$.
These are precisely the conditions used later to count the admissible
values of $b_a$.

Put $K=\F_{p^h}$ and consider the subrings
\[
 S=K[y;\sigma^d],\quad y=x^d,
 \qquad T=K[t;\sigma^n],\quad t=x^n=y^m.
\]
Both automorphisms on $K$ have order $H$ and fixed field $\F_P$.
The corresponding central variables satisfy
\[
 Z=x^h=y^H,\qquad W=t^H=Z^m.
\]
The inclusions and the corresponding centers can be summarized as
\[
\begingroup
\renewcommand{\arraystretch}{1.3}
\begin{array}{c|c|c}
 \text{ring}&\text{center}&\text{central variable}\\ \hline
 K[x;\sigma]&\F_p[Z]&Z=x^h\\
 S=K[y;\sigma^d]&\F_P[Z]&Z=y^H\\
 T=K[t;\sigma^n]&\F_P[W]&W=t^H=Z^m
\end{array}
\endgroup
\]
The passage from $T$ to $S$ substitutes $W=Z^m$; the passage from
$S$ to the original ring takes the product of the $d$ coefficient
conjugates over $\F_p$. The next lemma makes these two operations precise.

For $a\in K$, define $b_a\in\F_P$ by
\begin{equation}\label{eq:quadraticnorm}
 \Nrd_T(t^2+at+1)=G_a(W)=W^2+b_aW+1.
\end{equation}
This norm has the stated shape: the companion matrix is
\[
 A_a=\begin{pmatrix}0&-1\\1&-a\end{pmatrix},
\]
and the $H$-fold iterate of the associated semilinear map is
$K$-linear with determinant one. Its characteristic polynomial
therefore has degree two and constant term one.

\begin{lemma}\label{lem:normidentity}
In the original skew polynomial ring,
\[
 \Nrd_{K[x;\sigma]}(x^{2n}+ax^n+1)
 =\prod_{i=0}^{d-1}(Z^{2m}+b_a^{p^i}Z^m+1).
\]
\end{lemma}
\begin{proof}
\noindent\emph{Step 1: compare the norms in $T$ and $S$.}
We claim that
\[
 \Nrd_S(x^{2n}+ax^n+1)=G_a(Z^m).
\]
Indeed, extend the center of the cyclic algebra obtained by centrally
localizing $T$ from $\F_P(W)$ to $\F_P(Z)$ by $W\mapsto Z^m$.
The map fixes $K$ and sends $t$ to $y^m$; it respects the skew
relation because $y^m a=\sigma^n(a)y^m$ for $a\in K$.
It identifies the resulting cyclic algebra with that of $S$.
Indeed, since $\gcd(m,H)=1$, choose integers $u,v$ with $um+vH=1$.
After localization, $t$ and $Z$ are invertible, since $t^H=W$ and
$W,Z$ are nonzero elements of the respective centers' fraction fields.
Thus negative exponents are allowed, and the image contains
$(y^m)^uZ^v=y$. The map is therefore surjective.
Both algebras have dimension $H^2$ over $\F_P(Z)$, so the map is an
isomorphism. Reduced norms commute with extension of the center,
including inseparable extensions. This proves the first assertion as
an identity of rational functions and hence of polynomials.

\smallskip
\noindent\emph{Step 2: recover the norm in the original ring.}
Compute this reduced norm using right multiplication in
the left $K(Z)$-basis $\{1,x,\ldots,x^{h-1}\}$.
Group the basis vectors according to their exponents modulo $d$.
Since all exponents of $x^{2n}+ax^n+1$ are multiples of $d$, the matrix
is block diagonal, with $d$ blocks of size $H$.
The block indexed by $i$ is obtained from the matrix over $S$ by
applying $\sigma^i$ to its coefficients. Its determinant is
$G_a^{\sigma^i}(Z^m)$. Multiplying the block determinants proves the claim.
\end{proof}

\begin{lemma}\label{lem:composition}
For $b\in\F_{p^d}$ put
\[
 J_b(V)=\prod_{i=0}^{d-1}(V^2+b^{p^i}V+1)\in\F_p[V].
\]
Then $J_b(Z^m)$ is irreducible in 
$\F_p[Z]$ if and only if
\begin{enumerate}
\item $V^2+bV+1$ is irreducible over $\F_{p^d}$;
\item $[\F_p(b):\F_p]=d$;
\item for a root $\alpha$ of this quadratic, $Z^m-\alpha$ is irreducible
      over $\F_{p^{2d}}$.
\end{enumerate}
\end{lemma}
\begin{proof}
\noindent\emph{Step 1: determine when $J_b$ is irreducible.}
If $b$ belongs to a proper subfield $\F_{p^e}$, $e\mid d$, then $J_b$
is a $d/e$-th power. If the quadratic splits over $\F_{p^d}$, taking
coefficient norms of its two linear factors gives a factorization of
$J_b$ into two polynomials of degree $d$. Thus irreducibility of $J_b$
requires the first two conditions.

Conversely, under those two conditions, $b=-(\alpha+\alpha^{-1})$
shows that $\F_p(\alpha)$ contains $\F_{p^d}$. The quadratic is
irreducible there, so $[\F_p(\alpha):\F_p]=2d$. Hence $J_b$ is the
minimal polynomial of $\alpha$.

\smallskip
\noindent\emph{Step 2: account for the substitution $V=Z^m$.}
If $J_b$ is reducible, then substituting $Z^m$ in its factors shows
that $J_b(Z^m)$ is reducible. Thus it remains to work under the first
two conditions. Choose $\beta$ with $\beta^m=\alpha$. Since
$\alpha\in\F_p(\beta)$, the tower formula gives
\[
 [\F_p(\beta):\F_p]
 =2d\,[\F_p(\beta):\F_{p^{2d}}].
\]
The polynomial $J_b(Z^m)$ has degree $2dm$ and has $\beta$ as a root.
It is therefore irreducible exactly when the left-hand side is
$2dm$, or equivalently when $\beta$ has degree $m$ over
$\F_{p^{2d}}$. This is exactly the irreducibility of $Z^m-\alpha$.
\end{proof}

\section{Uniformity for the standard Frobenius}\label{sec:standard}
For this section let $P$ be any prime power, let $K=\F_{P^k}$, and
use the automorphism $c\mapsto c^P$. 
Put
\[
 T_1=K[t;\tau_1],\qquad \tau_1(c)=c^P,\qquad Z=t^k.
\]
Define $b_1(a) \in \F_P$ by
\[
 \Nrd_{T_1}(t^2+at+1)=Z^2+b_1(a)Z+1,
 \qquad w_{k,1}(b)=\#\{a\in K:b_1(a)=b\}.
\]

\begin{proposition}\label{prop:standard}
If $Z^2+bZ+1$ is irreducible over $\F_P$, then
$w_{k,1}(b)=D_k(P)$.
\end{proposition}
\begin{proof}
We first realize the desired fiber count as a fixed-point count on a
curve: equation~\eqref{eq:fixedfiber} identifies $(P+1)w_{k,1}(b)$
with the number of fixed points of a twisted Frobenius on the Drinfeld
curve $Y$. 
We then pass to its smooth projective closure and apply the Grothendieck--Lefschetz fixed-point formula. The key simplification is that the
twisted square of Frobenius has only the eigenvalue $-P$ on first
cohomology. Consequently, all the required traces reduce to the two
initial cases of an elliptic group element and such an element
composed with Frobenius. This gives the explicit value $D_k(P)$.

\smallskip
\noindent\emph{Step 1: realize the norm fiber as a fixed-point set.}
Let $F$ denote the coordinatewise $P$-power Frobenius over
$\overline{\F}_P$, and consider the affine Drinfeld curve
\[
 Y=\left\{u=(u_1,u_2)^t:\det(u,Fu)=1\right\}
\]
over $\F_P$, which 
admits a natural $\SL_2(\F_P)$-action. 
Define the morphism
\[
 \Gamma:Y\longrightarrow\mathbb A^1_{\F_P},
 \qquad
 u\longmapsto-\det(u,F^2u).
\]
Since
\[
 \det(Fu,F^2u)=\det(u,Fu)^P=1,
\]
expressing $F^2u$ in the basis $(u,Fu)$ gives
\begin{equation}\label{eq:linearized}
(F^2+\Gamma(u)F+1)u=0.
\end{equation}
For $a\in\overline{\F}_P$, we write $\Gamma^{-1}(a)$ for the geometric fiber in $Y(\overline{\F}_P)$.

Each geometric fiber of $\Gamma$ is a simply transitive $\SL_2(\F_P)$-set under the action on column vectors.
Indeed, let $V_a$ be the root space of
\[
 X^{P^2}+aX^P+X.
\]
This polynomial is separable and $\F_P$-linear, so $V_a$
is a two-dimensional $\F_P$-vector space.
For an ordered basis $(u_1,u_2)$ of $V_a$, its Moore determinant
\[
 D=u_1u_2^P-u_1^Pu_2
\]
is nonzero. The root relation gives $D^P=D$, so
$D\in\F_P^\times$. Rescaling one basis vector therefore
produces a basis with $D=1$.
By \eqref{eq:linearized}, the points of $\Gamma^{-1}(a)$
are precisely the column vectors $(u_1,u_2)^t$ arising
from these normalized ordered bases of $V_a$.
Any two such bases are related by a unique element
of $\SL_2(\F_P)$, proving the assertion.

We will also use the identities
\[
 \Gamma(F^ku)=\Gamma(u)^{P^k},
 \qquad
 \Gamma(gu)=\Gamma(u)
 \quad(g\in\SL_2(\F_P)).
\]
The first follows because $\Gamma$ is defined over $\F_P$.
For the second, $F^2(gu)=gF^2u$ and $\det g=1$ give
\[
 \Gamma(gu)
 =-\det(gu,gF^2u)
 =\Gamma(u).
\]

Choose $C\in\SL_2(\F_P)$ with characteristic polynomial
\[
 G_b(Z)=Z^2+bZ+1.
\]
For $a\in K=\F_{P^k}$ and $u\in\Gamma^{-1}(a)$, we have
\[
 \Gamma(F^ku)=a^{P^k}=a.
\]
Thus $F^ku$ belongs to the same fiber as $u$.
By simple transitivity, there is a unique
$C_u\in\SL_2(\F_P)$ such that
\[
 F^ku=C_u u.
\]
Put
\[
 B=[u,Fu],
 \qquad
 A_a=
 \begin{pmatrix}
 0&-1\\
 1&-a
 \end{pmatrix}.
\]
Equation \eqref{eq:linearized} gives
\[
 F(B)=BA_a.
\]
Since $F$ commutes with $C_u$, iteration yields
\[
 F^k(B)
 =B A_aF(A_a)\cdots F^{k-1}(A_a)
 =C_u B.
\]
The matrix $B$ is invertible because $\det(u,Fu)=1$.
Consequently, $C_u$ is conjugate to
$A_aF(A_a)\cdots F^{k-1}(A_a)$.
By the companion description of the reduced norm, its
characteristic polynomial is
\[
 Z^2+b_1(a)Z+1.
\]

Matrices in $\SL_2(\F_P)$ with irreducible characteristic
polynomial $G_b$ form one conjugacy class.
Their centralizer in $\SL_2(\F_P)$ is isomorphic to the
norm-one subgroup of $\F_{P^2}^{\times}$ and therefore
has order $P+1$.
For $g\in\SL_2(\F_P)$, the identity
\[
 F^k(gu)=gC_u u=(gC_u g^{-1})(gu)
\]
shows that
\[
 C_{gu}=gC_u g^{-1}.
\]
It follows that, in a fiber with $b_1(a)=b$, exactly $P+1$
points $u$ satisfy $C_u=C$, or equivalently
$C^{-1}F^ku=u$.
Fibers with $b_1(a)\ne b$ contribute no such points.

Conversely, if $u\in Y(\overline{\F}_P)$ satisfies
$C^{-1}F^ku=u$, then
\[
 \Gamma(u)^{P^k}
 =\Gamma(F^ku)
 =\Gamma(Cu)
 =\Gamma(u).
\]
Hence $\Gamma(u)\in K$, so every fixed point belongs to
one of the fibers considered above.
Therefore
\begin{equation}\label{eq:fixedfiber}
 (P+1)w_{k,1}(b)=\#Y^{C^{-1}F^k},
\end{equation}
where the right-hand side denotes the number of geometric
fixed points of $C^{-1}F^k$ on $Y$.
In particular, for $k=1$ this number is $P+1$, since
$b_1(a)=a$ and hence $w_{1,1}(b)=1$.

\smallskip
\noindent\emph{Step 2: compute the action of a twisted square of Frobenius.}
The smooth projective closure is
\[
 \overline Y:\quad S T^P- S^P T=U^{P+1}.
\]
The affine curve $Y$ is identified with the chart $U\ne0$
via $u_1=S/U$ and $u_2=T/U$.
The partial derivatives of its defining equation are
$T^P,-S^P,-U^P$, which have no common projective zero.
It is thus a smooth plane curve of degree $P+1$, with genus
$g=P(P-1)/2$. Setting $U=0$ identifies its boundary with
$\mathbb P^1(\F_P)$.
Write $\varepsilon=-I_2$ and
$\mathcal H=H^1_{\mathrm{\acute et}}
(\overline{Y}_{\overline{\F}_P},\overline{\mathbb Q}_\ell)$,
where $\ell$ is a prime different from the characteristic.
Here $F$ is the curve morphism raising coordinates to the $P$-th
power, and its action on cohomology is the pullback $F^*$, with
Weil eigenvalues of absolute value $P^{1/2}$ on $\mathcal H$.
We use the same notation for commuting curve endomorphisms and
their pullback actions on cohomology.
For $h\in\SL_2(\F_P)$ and $j\ge1$, the map $hF^j$ has zero
differential, so its fixed points are simple.
Its traces on $H^0$ and $H^2$ are $1$ and $P^j$, respectively.
Thus the trace formula reads
\[
 \#\overline Y^{hF^j}
 =1+P^j-\Tr(hF^j\mid\mathcal H).
\]

There are $P^3-P$ fixed points of $\varepsilon F^2$ on $Y$.
Indeed, $u_1^{P^2}=-u_1$ has $P^2-1$ nonzero solutions.
For each such $u_1$, the equation
$u_1u_2^P-u_1^Pu_2=1$ is a degree-$P$ polynomial equation in
$u_2$ with nonzero derivative $-u_1^P$, so it has exactly $P$
distinct solutions over $\overline{\F}_P$.
Raising this equation to the $P$-th power and using
$u_1^{P^2}=-u_1$ gives
$u_1^Pu_2^{P^2}+u_1u_2^P=1$.
Subtracting the original equation yields
$u_1^P(u_2^{P^2}+u_2)=0$, hence $u_2^{P^2}=-u_2$.
No solution has $u_1=0$. 
Adding the $P+1$ boundary points gives $P^3+1$.
The Grothendieck--Lefschetz trace formula gives
\[
 \Tr(\varepsilon F^2\mid\mathcal H)
 =1+P^2-(P^3+1)
 =-P^2(P-1)
 =-2gP.
\]
To determine the absolute values of its eigenvalues, put
$\Psi=\varepsilon F^2$. Since $\varepsilon$ commutes with $F$
and $\varepsilon^2=1$, we have
\[
 \Psi^2=F^4.
\]
By the Weil conjectures for smooth projective curves,
every eigenvalue of $F^4$ on $\mathcal H$ has complex
absolute value $P^2$.
Hence every eigenvalue of $\Psi$ has complex absolute value $P$.
There are $2g$ such eigenvalues, counted with multiplicity,
and their sum is $-2gP$. Equality in the triangle inequality
therefore forces every eigenvalue of $\Psi$ to be $-P$.
In particular, for every operator $A$ commuting with $\varepsilon F^2$,
\begin{equation}\label{eq:tracepower}
 \Tr\bigl((\varepsilon F^2)^sA\mid\mathcal H\bigr)
   =(-P)^s\Tr(A\mid\mathcal H).
\end{equation}
Indeed, after extending the coefficient field to an algebraic closure,
the commuting operators $\Psi$ and $A$ can be simultaneously upper
triangularized. Every diagonal entry of $\Psi$ is $-P$, so the diagonal
entries of $\Psi^sA$ are $(-P)^s$ times those of $A$.
No semisimplicity assumption is needed.

\smallskip
\noindent\emph{Step 3: obtain the two initial traces.}
If $g_0\in\SL_2(\F_P)$ has irreducible characteristic polynomial,
it has no fixed points on $\overline Y$. On the affine curve a fixed
vector would give eigenvalue one, and on the boundary a fixed point
would give an $\F_P$-rational eigenline.
Since $g_0$ has order prime to the characteristic and
has no fixed points on $\overline Y$, the
Grothendieck--Lefschetz trace formula gives
\[
 0=1-\Tr(g_0\mid\mathcal H)+1,
\]
and hence $\Tr(g_0\mid\mathcal H)=2$.
Equation \eqref{eq:fixedfiber} for $k=1$ gives
$\#Y^{g_0F}=P+1$. There are no fixed boundary points, so
\[
 \Tr(g_0F\mid\mathcal H)=0.
\]
Both statements remain true after multiplying $g_0$ by $\varepsilon$.

\smallskip
\noindent\emph{Step 4: evaluate the trace for every $k$.}
Put $\Psi=\varepsilon F^2$. Since $\varepsilon^2=1$ and all the
operators in question commute, we have
\[
 C^{-1}F^{2s}=\Psi^s(C^{-1}\varepsilon^s),\qquad
 C^{-1}F^{2s+1}=\Psi^s(C^{-1}\varepsilon^s F).
\]
The element $C^{-1}\varepsilon^s$ again has irreducible characteristic
polynomial. Hence Step 3 and \eqref{eq:tracepower} give trace
$2(-P)^s$ in the even case and zero in the odd case.
The action on the boundary is that of $C^{-1}$, since Frobenius fixes
$\mathbb P^1(\F_P)$ pointwise; it has no fixed points there.
Applying the Grothendieck--Lefschetz formula on $\overline Y$ therefore gives
\[
 \#Y^{C^{-1}F^k}=
 \begin{cases}
 P^k+1,&k\text{ odd},\\
 P^k+1-2(-P)^{k/2},&k\text{ even}.
 \end{cases}
\]
Division by $P+1$ in \eqref{eq:fixedfiber} proves the proposition.
\end{proof}

\section{Changing the Frobenius generator}\label{sec:generator}

Let $K=\F_{P^k}$ and let
\[
 \tau(c)=c^{P^r},\qquad \gcd(r,k)=1.
\]
For $k>1$, we may and do assume $1\le r<k$.
Put \(T_r=K[t;\tau]\) and \(Z=t^k\).
Define \(b_r(a)\in\F_P\) by
\[
 \Nrd_{T_r}(t^2+at+1)=Z^2+b_r(a)Z+1,
\]
and set
\[
 w_{k,r}(b)=\#\{a\in K:b_r(a)=b\}.
\]

The purpose of this section is to prove that changing the generator
of $\Gal(K/\F_P)$ does not change the fiber sizes over irreducible
quadratics. More precisely,
we establish the following reduction to Section~\ref{sec:standard}.

\begin{proposition}\label{prop:allgenerators}
For any $r$ relatively prime to $k$,
and any $b\in\F_P$ for which
$Z^2+bZ+1$ is irreducible over $\F_P$,
\[
 w_{k,r}(b)=w_{k,1}(b)=D_k(P).
\]
\end{proposition}

For $k=1$ this is already Proposition~\ref{prop:standard}; henceforth assume
$k>1$. We first establish the conjugacy classification and normal forms,
then express a fiber size as a count of normalized cyclic vectors.
Finally, we compare that count with the standard-Frobenius count,
using power maps with the same image and the same number of preimages.

\subsection{Conjugacy and normal forms}\label{subsec:classification}

In standard coordinates, a $\tau$-semilinear map has the form
$F(v)=A\tau(v)$, where $\tau$ acts coordinatewise.
We call $A$ its coefficient matrix. Its iterates satisfy
\[
 F^i(v)=A\tau(A)\cdots\tau^{i-1}(A)\tau^i(v).
\]
In particular, $F^k$ is $K$-linear because $\tau^k=1$.
A $\tau$-semilinear map $F$ on $K^2$ is said to have
determinant one if
\[
 \omega(Fv,Fw)=\tau(\omega(v,w)),
 \qquad \omega(v,w)=\det(v,w).
\]
Equivalently, the coefficient matrix of $F$ has
determinant one.

For $b\in\F_P$ such that
\[
 G_b(Z)=Z^2+bZ+1
\]
is irreducible, let $\mathcal F_{k,r}(b)$ denote the
set of determinant-one $\tau$-semilinear maps $F$
on $K^2$ for which $F^k$ has characteristic
polynomial $G_b$.
Here $\SL_2(K)$ acts on semilinear maps by conjugation:
two maps $F$ and $F'$ are conjugate if
\[
 F'=gFg^{-1}
 \qquad\text{for some }g\in\SL_2(K).
\]
Equivalently, if $F=A\tau$ and $F'=A'\tau$, then
\[
 A'=gA\tau(g)^{-1},
\]
where $\tau(g)$ is obtained by applying $\tau$ to each entry of $g$.

The conjugacy question is a twisted conjugacy problem: it belongs
to the framework of Lang--Steinberg theory and Shintani descent;
see, for example, \cite[Section~2.2]{Harper} for the correspondence
of classes and centralizers. We give an explicit two-dimensional
proof to obtain the normal forms used in the counting argument.
The classification and centralizer calculation are structural
ingredients; the additional task is to count normalized cyclic
vectors, or equivalently intersections with the companion family
$\{A_a:a\in K\}$, and compare these counts for different generators.

\begin{lemma}\label{lem:classes}
The set $\mathcal F_{k,r}(b)$ consists of a single
$\SL_2(K)$-conjugacy class. The centralizer in
$\SL_2(K)$ of any $F\in\mathcal F_{k,r}(b)$
has order $P+1$.

\end{lemma}

\begin{proof}
We prove the conjugacy assertion in three steps.
First, we derive normal forms and use them to show that
$\mathcal F_{k,r}(b)$ is nonempty.
Second, we show that any two such maps are
$\GL_2(K)$-conjugate and compute their invertible
self-intertwiners. Here an invertible self-intertwiner of $F$
means an element $g\in\GL_2(K)$ satisfying $gF=Fg$.
Finally, we use the determinants of these self-intertwiners
to obtain $\SL_2(K)$-conjugacy and determine the centralizer
in $\SL_2(K)$.

Let $\alpha$ be a root of $G_b$ in $\F_{P^2}$. Then
\[
 \alpha^P=\alpha^{-1},
 \qquad \alpha^{P+1}=1,
\]
and $\alpha\notin\F_P$.

We classify the semilinear maps in $\mathcal F_{k,r}(b)$.

\smallskip
\noindent
\emph{Case 1: $k=2s$ is even.}

Since $k$ is even and $\gcd(r,k)=1$, the integer $r$
is odd. Thus $\tau$ acts on $\F_{P^2}$ as the
$P$-power Frobenius, and in particular
\[
 \tau(\alpha)=\alpha^P=\alpha^{-1}.
\]
The two eigenvalues of $F^k$ belong to $K$ and are distinct.
Order the eigenlines so that the first eigenvalue is $\alpha$.
The identity $F^kF=FF^k$ and the $\tau$-semilinearity of $F$
show that $F$ exchanges these eigenlines: if $F^kv=\alpha v$, then
\[
 F^k(Fv)=F(F^kv)=\tau(\alpha)Fv=\alpha^{-1}Fv.
\]
Choose eigenvectors $e_1,e_2$ with
\[
 F^ke_1=\alpha e_1,\qquad F^ke_2=\alpha^{-1}e_2,
 \qquad \det(e_1,e_2)=1.
\]
Such a normalization is possible by rescaling one eigenvector.
Since $F$ is invertible and exchanges the eigenlines, write
\[
 F(e_1)=d e_2,\qquad F(e_2)=c e_1,
 \qquad c,d\in K^\times.
\]
Semilinearity now gives
\[
 F(ue_1+ve_2)=cv^\tau e_1+du^\tau e_2.
\]
In this determinant-one basis, the coefficient matrix of $F$ is
$\left(\begin{smallmatrix}0&c\\d&0\end{smallmatrix}\right)$.
Its determinant is $-cd=1$, so $d=-c^{-1}$.
Thus, in the resulting coordinates,
\begin{equation}\label{eq:evenform}
 F(u,v)=(cv^\tau,-c^{-1}u^\tau),
 \qquad c\in K^\times.
\end{equation}

Put
\[
 N=\Nm_{K/\F_{P^2}}(c).
\]
A direct iteration gives
\[
 F^2(u,v)=
 \left(-c\,\tau(c)^{-1}u^{\tau^2},
       -c^{-1}\tau(c)v^{\tau^2}\right).
\]
Since $\tau^2$ generates
$\Gal(K/\F_{P^2})$, iterating this identity
$s$ times yields
\[
 F^k=
 \operatorname{diag}
 \bigl((-1)^sN^{1-P},
       (-1)^sN^{P-1}\bigr).
\]
Thus $F^k$ has characteristic polynomial
$G_b$ precisely when its eigenvalues are
$\alpha$ and $\alpha^{-1}$.

The norm map $K^\times\to\F_{P^2}^\times$
is surjective, and the map
\[
 N\longmapsto N^{1-P}
\]
has image $\mu_{P+1}$, the group of $(P+1)$-st roots of unity
in $\F_{P^2}^{\times}$.
Consequently, we may choose $N$ so that
\[
 N^{1-P}=(-1)^s\alpha.
\]
This proves existence.

Next consider two maps of the form
\eqref{eq:evenform}, with parameters $c,c'$,
whose first eigenvalues agree. Set
\[
 \delta=c'/c.
\]
The equality of their first eigenvalues gives
\[
 \Nm_{K/\F_{P^2}}(\delta)^{P-1}=1,
\]
hence
\[
 \Nm_{K/\F_{P^2}}(\delta)\in\F_P^\times.
\]
It follows that
\[
 \Nm_{K/\F_{P^2}}
 \left(\frac{\tau(\delta)}{\delta}\right)=1.
\]
By Hilbert's Theorem 90 for the cyclic
extension $K/\F_{P^2}$, whose Galois
group is generated by $\tau^2$, there is
$u\in K^\times$ such that
\[
 \frac{\tau^2(u)}u
 =\frac{\tau(\delta)}{\delta}.
\]

Define
\[
 g=\operatorname{diag}
       \bigl(u,\delta^{-1}\tau(u)\bigr).
\]
The defining equation for $u$ gives $gF_c=F_{c'}g$.
The two maps are therefore $\GL_2(K)$-conjugate.

Any invertible self-intertwiner $g$ of $F_c$ also commutes with
$F_c^k$. Since $F_c^k$ is diagonal with distinct eigenvalues, $g$
must be diagonal. Writing $g=\operatorname{diag}(u,v)$, the identity
$gF_c=F_cg$ is equivalent to
\[
 u=\tau(v),\qquad v=\tau(u).
\]
Thus $\tau^2(u)=u$. Since the fixed field of $\tau^2$ is
$\F_{P^2}$, the invertible self-intertwiners are exactly
\[
 \operatorname{diag}(u,\tau(u)),
 \qquad u\in\F_{P^2}^\times.
\]
Their determinants are
\[
 u\tau(u)=u^{P+1}.
\]
Hence the determinant map on the
self-intertwiner group is the surjective norm
\[
 \Nm_{\F_{P^2}/\F_P}:
 \F_{P^2}^\times\longrightarrow\F_P^\times.
\]

\smallskip
\noindent
\emph{Case 2: $k$ is odd.}

Put
\[
 Q=P^k,\qquad E=\F_{P^{2k}}.
\]
Since $k$ is odd, $K=\F_Q$ does not contain
$\F_{P^2}$. Therefore $G_b$ remains
irreducible over $K$.

We may identify the two-dimensional $K$-vector space $K^2$ with $E$ so that
$F^k$ acts as multiplication by $\alpha$.

Let $j$ be the unique even integer satisfying
\[
 0\le j<2k,\qquad j\equiv r\pmod{k},
\]
and put
\[
 \rho(z)=z^{P^j}.
\]
Then $\rho|_K=\tau$, and $\rho$ fixes
$\F_{P^2}$ pointwise. Moreover,
\[
 \gcd(j/2,k)=1,
\]
so $\rho$ generates
$\Gal(E/\F_{P^2})$.

Because $F$ commutes with $F^k$, it commutes
with multiplication by $\alpha$.
As $E=K(\alpha)$, the map $F$ is
$\rho$-semilinear on the one-dimensional
$E$-vector space $E$. Consequently,
\begin{equation}\label{eq:oddform}
 F(z)=c\rho(z)
\end{equation}
for some $c\in E^\times$.

Since $\rho$ generates
$\Gal(E/\F_{P^2})$, iteration gives
\[
 F^k(z)=
 \Nm_{E/\F_{P^2}}(c)\,z.
\]
Thus the condition $F^k(z)=\alpha z$ for every $z\in E$ is
equivalent to
\[
 \Nm_{E/\F_{P^2}}(c)=\alpha.
\]

We next express the determinant-one condition in these
coordinates, before verifying that both conditions on $c$ can be satisfied simultaneously.
Choose $\theta\in\F_{P^2}^\times$ such that
$\theta^P=-\theta$; in characteristic two, take $\theta=1$.
Define
\[
 \omega_\theta(z,w)=\theta(zw^Q-z^Qw).
\]
This is a nondegenerate, alternating,
$K$-valued bilinear form on $E$.

Indeed, $\theta^Q=-\theta$, since $k$ is odd, and
therefore
\[
 \omega_\theta(z,w)^Q=\omega_\theta(z,w).
\]

Let $\phi:E\longrightarrow K^2$ be an identification for which
$\phi^{-1}F^k\phi$ is multiplication by $\alpha$.
The space of alternating $K$-bilinear forms on the two-dimensional
space $E$ is one-dimensional. Since both forms are nondegenerate,
there is $\lambda\in K^\times$ such that
\[
 \det(\phi(z),\phi(w))=\lambda\omega_\theta(z,w).
\]
Multiplication by $u\in E^\times$ satisfies
\[
 \omega_\theta(uz,uw)=u^{Q+1}\omega_\theta(z,w)
              =\Nm_{E/K}(u)\omega_\theta(z,w).
\]
Choose $u$ with $\Nm_{E/K}(u)=\lambda^{-1}$, which is possible
by surjectivity of the norm, and replace $\phi(z)$ by $\phi(uz)$.
The transported determinant form is then exactly $\omega_\theta$.
This change preserves the description of $F^k$ as multiplication
by $\alpha$, since multiplication by $u$ commutes with
multiplication by $\alpha$.
In these normalized coordinates, write $F=c\rho$ again.

As $\rho(\theta)=\theta$ and $\rho$ commutes with
the $Q$-power Frobenius, we have
\[
 \omega_\theta(Fz,Fw)
 =c^{Q+1}\rho(\omega_\theta(z,w)).
\]
Thus $F$ has determinant one if and only if
\[
 \Nm_{E/K}(c)=c^{Q+1}=1.
\]

We now verify that both norm conditions can be satisfied simultaneously.

Let
\[
 U=\ker(\Nm_{E/K}).
\]
This is a cyclic group of order $Q+1$.
Set
\[
 A=\frac{Q+1}{P+1},
 \qquad
 B=\frac{Q-1}{P-1}.
\]
The exponent of $\Nm_{E/\F_{P^2}}$
is
\[
 \frac{P^{2k}-1}{P^2-1}=AB.
\]
Since $k$ is odd,
\[
 B=1+P+\cdots+P^{k-1}
 \equiv1\pmod{P+1}.
\]
Hence
\[
 \gcd(B,P+1)=1,
 \qquad
 \gcd(Q+1,AB)=A.
\]
It follows that the image of the norm map
on $U$ has order
\[
 \frac{Q+1}{A}=P+1.
\]
This image is contained in
$\mu_{P+1}\subset\F_{P^2}^\times$.
Therefore
\[
 \Nm_{E/\F_{P^2}}(U)=\mu_{P+1}.
\]
In particular, there exists $c\in U$ with
\[
 \Nm_{E/\F_{P^2}}(c)=\alpha.
\]
This proves existence.

Suppose now that $F_c=c\rho$ and
$F_{c'}=c'\rho$ have the same $k$-th power.
Then
\[
 \Nm_{E/\F_{P^2}}(c'/c)=1.
\]
Hilbert's Theorem 90 gives $u\in E^\times$
such that
\[
 \frac{\rho(u)}u=\frac c{c'}.
\]
Multiplication by $u$ intertwines $F_c$
and $F_{c'}$, proving $\GL_2(K)$-conjugacy.

Any self-intertwiner commutes with multiplication by $\alpha$.
Since $E=K(\alpha)$, it is therefore $E$-linear and hence is
multiplication by some $u\in E$. Commutation with $c\rho$ then
requires $\rho(u)=u$. Consequently, the invertible
self-intertwiners are exactly the multiplications by
\[
 u\in E^{\rho,\times}=\F_{P^2}^\times.
\]
Their determinants as $K$-linear maps are
\[
 \Nm_{E/K}(u)=u^{Q+1}=u^{P+1},
\]
where the last equality follows because
$k$ is odd.

Thus, also in this case, the determinant
map on the self-intertwiner group is the
surjective norm
\[
 \Nm_{\F_{P^2}/\F_P}.
\]

\smallskip
\noindent
\emph{Passage to $\SL_2(K)$.}

In either case, let $g\in\GL_2(K)$
intertwine two determinant-one semilinear
maps. Taking determinants in the
intertwining identity gives
\[
 \tau(\det g)=\det g.
\]
Since the fixed field of $\tau$ is $\F_P$,
we have
\[
 \det g\in\F_P^\times.
\]

More explicitly, let $F,F'\in\mathcal F_{k,r}(b)$ and choose
$g\in\GL_2(K)$ with $gF=F'g$, as established above.
Put $\delta=\det g\in\F_P^\times$.
In either parity, the determinant map on the self-intertwiner group
of $F$ is surjective onto $\F_P^\times$.
We may therefore choose a $K$-linear self-intertwiner $u$ of $F$
with $\det u=\delta^{-1}$. Then
\[
 (gu)F=g(Fu)=F'(gu),\qquad \det(gu)=1.
\]
Thus any two elements of $\mathcal F_{k,r}(b)$ are
$\SL_2(K)$-conjugate. Since the constructions above also prove
that this set is nonempty, it is a single $\SL_2(K)$-conjugacy class.

In either parity, the self-intertwiner
group is isomorphic to
$\F_{P^2}^\times$, and its determinant
map is the norm to $\F_P^\times$.
Consequently, its intersection with
$\SL_2(K)$ has order
\[
 \frac{P^2-1}{P-1}=P+1.
\]
\end{proof}

\subsection{From norm fibers to cyclic vectors}

For $F\in\mathcal F_{k,r}(b)$, put
\[
 \mathcal C(F)=\{v\in K^2:\omega(v,Fv)=1\}.
\]
We claim that
\begin{equation}\label{eq:cyclic}
 w_{k,r}(b)=\frac{\#\mathcal C(F)}{P+1}.
\end{equation}

\begin{proof}
Let $F\in\mathcal F_{k,r}(b)$.
If $v\in\mathcal C(F)$, then $(v,Fv)$
is a determinant-one basis of $K^2$.

Since $F$ is $\tau$-semilinear and has
determinant one, there is a unique
$a\in K$ such that
\[
 F^2v=-v-aFv.
\]
Indeed, write $F^2v=\lambda v+\mu Fv$.
The determinant-one condition gives
\[
 1=\omega(Fv,F^2v)=-\lambda,
\]
so $\lambda=-1$; setting $a=-\mu$ gives the assertion.
Thus, relative to the basis $(v,Fv)$,
the map $F$ has coefficient matrix
\[
 A_a=
 \begin{pmatrix}
 0&-1\\
 1&-a
 \end{pmatrix}.
\]
The $K$-linear map $F^k$ has characteristic
polynomial
\[
 Z^2+b_r(a)Z+1.
\]
Since this polynomial is $G_b$, we obtain
$b_r(a)=b$.

Conversely, if $b_r(a)=b$, the semilinear map $A_a\tau$ belongs to $\mathcal
F_{k,r}(b)$ and hence is $\SL_2(K)$-conjugate to $F$.

The determinant-one bases that identify $F$ with $A_a\tau$ form a
torsor under the centralizer of $F$ in $\SL_2(K)$.
More explicitly, if $B_0$ is one such basis matrix, then all the
others are uniquely of the form $gB_0$ with
$g\in Z_{\SL_2(K)}(F)$, where $Z_{\SL_2(K)}(F)$ denotes this
centralizer. There are therefore exactly $P+1$ such bases.

Each such basis is uniquely determined by its first vector $v$, since the
second vector is $Fv$.

It follows that
\[
 \#\mathcal C(F)=(P+1)w_{k,r}(b),
\]
which proves \eqref{eq:cyclic}.

\end{proof}

\subsection{Comparison with the standard Frobenius}

\begin{proof}[Proof of Proposition~\ref{prop:allgenerators}]
By Lemma~\ref{lem:classes}, it suffices to compare the cyclic-vector
counts for $F$ and a standard-Frobenius map $F_0$ such that both
$F^k$ and $F_0^k$ have characteristic polynomial $G_b$.

\smallskip
\noindent
\emph{Case 1: $k=2s$ is even.}

Use the normal form \eqref{eq:evenform}:
\[
 F(u,v)=(cv^{P^r},-c^{-1}u^{P^r}).
\]
Define the standard-Frobenius semilinear map
\[
 F_0(u,v)=(cv^P,-c^{-1}u^P).
\]

Since $r$ is odd, the action of $\tau$ on
$\F_{P^2}$ is the $P$-power Frobenius.
Consequently, both $F^k$ and $F_0^k$
have eigenvalues
\[
 (-1)^sN^{1-P},
 \qquad
 (-1)^sN^{P-1},
 \qquad
 N=\Nm_{K/\F_{P^2}}(c).
\]
Thus $F\in\mathcal F_{k,r}(b)$ and $F_0\in\mathcal F_{k,1}(b)$.

Their cyclic-vector equations are
respectively
\[
 -c^{-1}u^{P^r+1}
 -cv^{P^r+1}=1
\]
and
\[
 -c^{-1}u^{P+1}
 -cv^{P+1}=1.
\]

Since $r$ is odd and $\gcd(r,k)=1$,
the standard gcd identity for powers gives
\[
 \gcd(P^r+1,P^k-1)=P+1.
\]
The multiplicative group $K^\times$
is cyclic. Therefore the power maps
\[
 x\longmapsto x^{P^r+1},
 \qquad
 x\longmapsto x^{P+1}
\]
have the same kernel size, namely $P+1$,
and hence the same image.

Moreover, every nonzero element of their common image has exactly $P+1$
preimages under either map, while zero has exactly one preimage.

It follows that the two cyclic-vector equations have the same number of
solutions.

\smallskip
\noindent
\emph{Case 2: $k$ is odd.}

We identify the two-dimensional $K$-vector space $K^2$ with the
quadratic extension $E/K$. We retain the notation
\[
 Q=P^k,\qquad E=\F_{P^{2k}},\qquad
 j\equiv r\pmod k,\quad j\ \text{even},\quad 0\le j<2k.
\]
Use the normal form \eqref{eq:oddform},
with the normalized alternating form
$\omega_\theta$:
\[
 F(z)=cz^{P^j},
 \qquad
 \Nm_{E/K}(c)=1,
 \qquad
 \Nm_{E/\F_{P^2}}(c)=\alpha.
\]

Define
\[
 \rho_0(z)=z^{P^{k+1}},
 \qquad
 F_0(z)=c\rho_0(z).
\]
Because $k$ is odd, $k+1$ is even and
\[
 \gcd(k+1,2k)=2.
\]
Thus $\rho_0$ generates
$\Gal(E/\F_{P^2})$, and its restriction
to $K$ is the standard $P$-power
Frobenius.

Both $F^k$ and $F_0^k$ are multiplication by $\alpha$. Furthermore, both maps
have determinant one with respect to $\omega_\theta$, because $\Nm_{E/K}(c)=1$.

Put
\[
 l=j+k.
\]
Since $j$ is even and $k$ is odd,
$l$ is odd. Also
\[
 \gcd(l,k)=\gcd(j,k)=1,
\]
so
\[
 \gcd(l,2k)=1.
\]

For $z\in E$, direct calculation gives
\[
\begin{aligned}
 \omega_\theta(z,Fz)
 &=\theta\left(
    c^Qz^{1+P^{j+k}}
    -cz^{P^k+P^j}
   \right)\\
 &=\theta(c^Qs-cs^Q),
\end{aligned}
\]
where
\[
 s=z^{P^l+1}.
\]
Here we used
\[
 (z^{P^l+1})^Q
 =z^{P^j+P^k},
\]
which holds because
$z^{P^{2k}}=z$.

For $F_0$, the same calculation gives
\[
 \omega_\theta(z,F_0z)
 =\theta(c^Qs_0-cs_0^Q),
 \qquad
 s_0=z^{P+1}.
\]

Since $l$ is odd and $\gcd(l,2k)=1$,
we have
\[
 \gcd(P^l+1,P^{2k}-1)=P+1.
\]
Consequently, the power maps
\[
 z\longmapsto z^{P^l+1},
 \qquad
 z\longmapsto z^{P+1}
\]
on $E^\times$ have the same image
and the same kernel size $P+1$.

As in the even case, zero has a unique
preimage under either power map.
Thus, for every $s\in E$, the numbers
of solutions of
\[
 z^{P^l+1}=s
 \quad\text{and}\quad
 z^{P+1}=s
\]
are equal.

Since the cyclic-vector conditions
for $F$ and $F_0$ are respectively
\[
 \theta(c^Qs-cs^Q)=1
\]
with $s=z^{P^l+1}$ and $s=z^{P+1}$,
their numbers of solutions coincide.

\smallskip
In both cases, the cyclic-vector count
for the generator $\tau$ agrees with
that for the standard Frobenius.

By Lemma~\ref{lem:classes},
\[
 w_{k,r}(b)=w_{k,1}(b).
\]
Proposition~\ref{prop:standard} gives
\[
 w_{k,1}(b)=D_k(P),
\]
and the assertion follows.
\end{proof}

\begin{remark}
The fiber counts satisfy $\sum_{b\in\F_P}w_{k,r}(b)=P^k$, since
each $a\in K$ has exactly one norm coefficient. There are
$(P-\epsilon_P)/2$ irreducible quadratics $Z^2+bZ+1$, where
$\epsilon_P$ is zero for even $P$ and one for odd $P$.
Consequently, the total contribution of irreducible norm values is
$(P-\epsilon_P)D_k(P)/2$, and that of reducible norm values is
$P^k-(P-\epsilon_P)D_k(P)/2$. The proposition concerns the individual
fibers over irreducible norm values.
\end{remark}

\section{Arithmetic count and conclusion}\label{sec:arithmetic}
We now count the admissible norm coefficients. These satisfy the
three conditions of Lemma~\ref{lem:composition}. We do this by counting roots
$\alpha$ in the norm-one group
\[
 \{\alpha\in\F_{P^2}^{\times}:\alpha^{P+1}=1\}.
\]
The binomial irreducibility criterion imposes conditions on the prime
parts of $\ord(\alpha)$, and M\"obius inversion selects the elements
of exact degree $2d$ over $\F_p$. Finally, we divide by two because
$\alpha$ and $\alpha^{-1}$ give the same coefficient
$b=-(\alpha+\alpha^{-1})$. This yields $B_{d,m}(p)$; multiplying by
the common fiber size from Section~\ref{sec:generator} completes the
proof of the main theorem.
\begin{lemma}\label{lem:arithmetic}
The number of $b\in\F_{p^d}$ satisfying the three conditions of
Lemma~\ref{lem:composition} is $B_{d,m}(p)$ of \eqref{eq:B}.
\end{lemma}
\begin{proof}
The first two conditions are equivalent to choosing
\[
 \alpha^{P+1}=1,\qquad [\F_p(\alpha):\F_p]=2d,
\]
and setting $b=-(\alpha+\alpha^{-1})$.
Each $b$ corresponds to exactly two such elements, $\alpha$ and
$\alpha^{-1}$. They are distinct, since an element equal to its inverse
has degree one over $\F_p$.

If $m$ is even, $\alpha$ is a square in $\F_{P^2}$:
this is automatic in characteristic two, while for odd $P$ one has
$\alpha^{(P^2-1)/2}=1$.
Writing $\alpha=\beta^2$, we obtain
\[
 X^m-\alpha
 =(X^{m/2}-\beta)(X^{m/2}+\beta),
\]
so $X^m-\alpha$ is reducible. This proves the assertion in this case.

Suppose $m>1$ is odd. The finite-field binomial criterion
\cite[Theorem 3.75]{LN} says that $X^m-\alpha$ is irreducible if and
only if every prime $\ell\mid m$ divides $\ord(\alpha)$ but does not
divide $(P^2-1)/\ord(\alpha)$; the additional condition for multiples
of four is irrelevant here.
Since $\ord(\alpha)\mid P+1$ and $\gcd(P-1,P+1)\mid2$,
this is exactly
\begin{equation}\label{eq:valuation}
 0<v_\ell(\ord(\alpha))=v_\ell(P+1)
 \qquad(\ell\mid m).
\end{equation}
For $m=1$, the binomial is linear and this condition is empty.
In particular, $R\nmid P+1$ implies that there are no such elements.

Assume $R\mid P+1$ and let $L$ be as in \eqref{eq:B}.
For each $e\mid2d$, the group
\[
 \left\{\alpha\in\F_{p^{2d}}^\times:\alpha^{P+1}=1\right\}
 \cap\F_{p^e}^{\times}
\]
is cyclic of order $g_e=\gcd(P+1,p^e-1)$.
It contains no elements satisfying \eqref{eq:valuation}
unless $L\mid g_e$.

Suppose $L\mid g_e$. For each prime $\ell\mid R$,
the $\ell$-primary component of this cyclic group has order
$\ell^{v_\ell(P+1)}$.
The proportion of elements whose orders have maximal
$\ell$-part is $1-1/\ell$.
The conditions for distinct primes are independent under
the decomposition into primary components.
Hence the number of elements satisfying \eqref{eq:valuation} is
\[
 g_e\prod_{\ell\mid R}(1-\ell^{-1})
 =g_e\frac{\varphi(R)}R.
\]
If $f=[\F_p(\alpha):\F_p]$, then
$\alpha\in\F_{p^e}$ if and only if $f\mid e$.
For $\alpha\in\F_{p^{2d}}$, the indicator that its degree over $\F_p$
is exactly $2d$ equals
\[
 \sum_{e\mid2d}\mu(2d/e)\,\mathbf1_{\alpha\in\F_{p^e}}.
\]
Summing this indicator over the elements satisfying
\eqref{eq:valuation}, and dividing by two for inversion, gives
\[
 B_{d,m}(p)
 =\frac{\varphi(R)}{2R}
   \sum_{\substack{e\mid2d\\L\mid g_e}}
   \mu(2d/e)g_e,
\]
as required.
\end{proof}

\begin{proof}[Proof of Theorem~\ref{thm:main}]
By Lemmas~\ref{lem:normidentity} and \ref{lem:composition}, the
admissible norm coefficients are precisely the $B_{d,m}(p)$ elements
counted in Lemma~\ref{lem:arithmetic}.
On $K=\F_{P^H}$ the automorphism $\sigma^n$ is
$c\mapsto c^{P^m}$, with $\gcd(m,H)=1$.
Proposition~\ref{prop:allgenerators}, applied with $k=H$ and $r=m$,
shows that each admissible coefficient has exactly $D_H(P)$ preimages.
Multiplication yields the asserted formula.
\end{proof}

\section{Specializations}\label{sec:specializations}
Recall that $\epsilon_p=0$ for even $p$ and $\epsilon_p=1$ for odd $p$.
For $m=1$ the arithmetic count simplifies to
\begin{equation}\label{eq:C}
 C_d(p):=B_{d,1}(p)
  =\frac12\sum_{\substack{e\mid d\\e\text{ odd}}}
            \mu(e)(p^{d/e}-\epsilon_p).
\end{equation}
One direct derivation is as follows. Put
\[
 U_j=\left\{c\in\F_{p^{2j}}^\times:c^{p^j+1}=1\right\}.
\]
The coefficients of irreducible quadratics $X^2+bX+1$ over
$\F_{p^j}$ are precisely the values $b=-(c+c^{-1})$ with
$c\in U_j$ and $c\ne c^{-1}$, each obtained from exactly two
such elements. The group $U_j$ has $p^j+1$ elements, of which
$1+\epsilon_p$ satisfy $c=c^{-1}$. Hence the number of these
quadratics is $(p^j-\epsilon_p)/2$.
An irreducible quadratic over a subfield remains irreducible
precisely under an extension of odd degree. Partitioning
coefficients by their exact degree gives
$(p^j-\epsilon_p)/2=\sum_{e\mid j,\ j/e\text{ odd}}C_e(p)$;
M\"obius inversion gives \eqref{eq:C}.

Let $S_p(2d)$ denote the number of monic self-reciprocal irreducible
polynomials of degree $2d$ over $\F_p$. Then
\begin{equation}\label{eq:selfreciprocal-count}
 C_d(p)=dS_p(2d).
\end{equation}
Indeed, the admissible coefficients counted by $C_d(p)$ correspond
to the inverse pairs $\{\alpha,\alpha^{-1}\}$ satisfying
\[
 \alpha^{p^d+1}=1,\qquad [\F_p(\alpha):\F_p]=2d,
\]
via $b=-(\alpha+\alpha^{-1})$.
The minimal polynomial of such an $\alpha$ is self-reciprocal,
since $\alpha^{-1}=\alpha^{p^d}$ is a Frobenius conjugate.
Conversely, if $\alpha$ is a root of a self-reciprocal irreducible
polynomial of degree $2d$, inversion permutes its Frobenius orbit
and is a nontrivial involution commuting with Frobenius.
It must therefore act on that orbit as the $d$-th power of Frobenius,
so $\alpha^{-1}=\alpha^{p^d}$.
Each such polynomial has $2d$ distinct roots, grouped into $d$
inverse pairs. Distinct inverse pairs give distinct values of $b$,
since their elements are precisely the roots of $X^2+bX+1$.
Thus each polynomial contributes exactly $d$ coefficients to $C_d(p)$,
proving \eqref{eq:selfreciprocal-count}.
Combining \eqref{eq:C} and \eqref{eq:selfreciprocal-count}
recovers the classical enumeration of monic self-reciprocal
irreducible polynomials; see \cite[Theorem 3]{MG}.

\begin{corollary}\label{cor:multiple}
For all $n,k\geq1$,
\[
 N_n(p,nk)=C_n(p)D_k(p^n)
         =nS_p(2n)D_k(p^n).
\]
If $n$ is a power of two, then
\[
 N_n(p,h)=
 \begin{cases}
 0,&n\nmid h,\\[2pt]
 \dfrac{p^n-\epsilon_p}{2}D_{h/n}(p^n),&n\mid h.
 \end{cases}
\]
In particular, this includes the counting formulas for $n=2$ and $n=4$.
\end{corollary}
\begin{proof}
If $h=nk$, then
\[
 d=\gcd(n,h)=n,\qquad m=1,\qquad H=k,\qquad P=p^n.
\]
Theorem~\ref{thm:main} therefore gives
\[
 N_n(p,nk)=D_k(p^n)B_{n,1}(p)
         =D_k(p^n)C_n(p).
\]
Equation~\eqref{eq:selfreciprocal-count} gives the expression in
terms of $S_p(2n)$.
For a power of two, $n\nmid h$ implies that
$m=n/\gcd(n,h)$ is even, so the count is zero by
Theorem~\ref{thm:main}. If $n\mid h$, the sum in
\eqref{eq:C} for $C_n(p)$ has only the odd divisor $e=1$,
which gives the stated formula.
\end{proof}

\begin{corollary}
For $n=3$,
\[
 N_3(p,h)=
 \begin{cases}
 0,&3\nmid h,\ p\equiv0,1\pmod3,\\[2pt]
 \dfrac{p+1}{3}D_h(p),&3\nmid h,\ p\equiv2\pmod3,\\[5pt]
 \dfrac{p^3-p}{2}D_{h/3}(p^3),&3\mid h.
 \end{cases}
\]
\end{corollary}
\begin{proof}
If $3\nmid h$, then $d=1,m=3$.
The arithmetic count vanishes unless $3\mid p+1$.
In that case the norm-one group has $p+1$ elements, and precisely
$2(p+1)/3$ have maximal $3$-part in their order. These elements have
degree two, since the elements of the norm-one group in $\F_p$ have
order at most two. Dividing by two gives $B_{1,3}(p)=(p+1)/3$.
If $3\mid h$, use $C_3(p)=(p^3-p)/2$.
\end{proof}

\begin{remark}
For an example with $n\nmid h$, take $(p,n,h)=(2,3,5)$.
Then $d=1$, $m=3$, $H=5$, and
\[
 B_{1,3}(2)=1,\qquad D_5(2)=\frac{2^5+1}{3}=11,
 \qquad N_3(2,5)=11.
\]
The automorphism in the quadratic reduction is $c\mapsto c^{2^3}$
on $\F_{2^5}$, a nonstandard Frobenius generator. This illustrates
the use of Section~\ref{sec:generator}.
\end{remark}

\section*{Declarations}

\subsection*{Competing interests}
The authors declare no competing interests.

\subsection*{Data availability}
No external datasets were used.

\subsection*{Declaration of generative AI and AI-assisted
technologies in the manuscript preparation process}

Before using generative AI in this work, the authors had proved
the counting formula for $n=2$ by a different method and had
formulated conjectures for the cases $n=3$ and $n=4$ based on
computations in Python.
These computations were performed on the computing cluster
at RIKEN iTHEMS.
They subsequently used OpenAI's Codex with the GPT-6.1 Sol model
to help formulate the conjectured counting formula for general $n$,
explore proof strategies, develop proof arguments, write code
for computational checks, and assist with drafting and revising
the manuscript.
The authors reviewed and edited all AI-assisted material and
take full responsibility for the mathematical content,
references, and final presentation.

\end{document}